\documentclass{amsart}
\usepackage{mathrsfs}
\usepackage{amsfonts}
\usepackage{amssymb}
\usepackage[papersize={8.3in,10.1in},textwidth=15.9cm,textheight=21.3cm,centering]{geometry}

\usepackage[colorlinks=true,citecolor=blue,linkcolor=blue]{hyperref}
\hypersetup{
pdfstartpage=1,
pdfstartview=FitH}

\allowdisplaybreaks

\newtheorem{theorem}{Theorem}
\newtheorem{lemma}{Lemma}

\newtheorem{definition}{Definition}

\theoremstyle{remark}

\newcommand{\ud}{\mathrm{d}}

\begin{document}

\title[Mean values of character sums analogue of Kloosterman sums]{Mean values of character sums analogue of\\ Kloosterman sums}
\thanks{The work is supported by N.S.F. (No. 11171265) of P.R. China.}
\author{Ping Xi}

\address{School of Mathematics and Statistics, Xi'an Jiaotong University, Xi'an 710049, P. R. China} \email{xprime@163.com}


\date{}


\begin{abstract}
Let $q$ be a positive integer, $\chi$ a nontrivial character mod
$q$, $\mathcal{I}$ an interval of length not exceeding $q.$ In this paper we shall study the
character sum analogue of the well-known Kloosterman sum,\[\sum_{\substack{a\in\mathcal{I}\\ \gcd(a,q)=1}}\chi(ma+n\overline{a}),\] where
$\overline{a}$ is the multiplicative inverse of $a\bmod q$. The mean square values and bilinear forms for such sums are proved.
\end{abstract}

\maketitle

\section{Introduction and statements of results}\label{sec:1}

Let $q$ be a positive integer, and $\chi$ a nontrivial
character mod $q$. A nontrivial bound of the character sum
\begin{equation}\label{eq:1}\sum_{x\in\mathcal{A}}\chi(f(x))\end{equation}can lead to important applications in number theory, where $f(x)$ is a rational function and the summation is over
a given set $\mathcal{A}$ avoiding the poles of $f$.

Estimates for character sums enjoy a fruitful literature due to their powerful applications. The classical result for $f$ being a linear function, independently proved by
G. P\'{o}lya \cite{P} and I. M. Vinogradov \cite{V}, is the
upper bound\[\sum_{n=N+1}^{N+H}\chi(n)\ll q^{1/2}\log q.\] Later, D. A. Burgess \cite{B1,B2} developed a powerful method for estimating the character sum, which resulted in the great improvement,
\begin{equation}\label{eq:2}\sum_{n=N+1}^{N+H}\chi(n)\ll
H^{1-1/r}q^{(r+1)/4r^2+\varepsilon}\end{equation}with $r=1,2,3$ for any $q$ and
with arbitrary positive integer $r$ if $q$ is cube-free. It should be mentioned that J. Friedlander and H. Iwaniec \cite{FI} have given a simple and illuminating proof for (\ref{eq:2}) using Fourier techniques and H\"{o}lder's inequality. Recently, A. Granville and K. Soundararajan \cite{GS} have managed to improve the P\'{o}lya-Vinogradov bound partially, obtaining
\[\sum_{n=N+1}^{N+H}\chi(n)\ll q^{1/2}(\log q)^{1-\delta_g/2+\varepsilon}\]for $\chi$ primitive mod $q$ and with odd order $g$, where $\delta_g=1-\frac{g}{\pi}\sin\frac{\pi}{g}$.

For the case of $f$ being an arbitrary polynomial with integral coefficients, one has the estimate
\[\sum_{x=N+1}^{N+H}\chi(f(x))\ll p^{1/2}\log p\]for $\chi$ being a
$d$-th order character mod $p$(prime) and $f(x)$ a polynomial
that is not a perfect $d$-th power (mod $p$), where the implied constant depends on the degree of $f(x)$. This is a well-known consequence of the theorem, due to A. Weil \cite{W},
that the Riemann Hypothesis is true for the zeta-function of an
algebraic function field over a finite field.

Another powerful analytic technique in number theory is the estimate for exponential sums, one of which is known as the Kloosterman sum, playing an important role in modern analytic number theory.
Let $c$ be a positive integer, for any integers $m,n$, the classical Kloosterman sum is defined by
\[\sum_{\substack{a\bmod c\\ \gcd(a,c)=1}}e\left(\frac{ma+n\overline{a}}{c}\right),\]where $e(x)=e^{2\pi
ix}$ and $\overline{a}$ is
the multiplicative inverse of $a\bmod c$.

In this paper, we consider a special and usually the most important case of (\ref{eq:1}), that is
\[f(x)=mx+\frac{n}{x}.\] To be precise, we shall study the character sum
\begin{equation}\label{eq:3}\Lambda_\chi(m,n,\mathcal{I};q)=\sum_{\substack{a\in\mathcal{I}\\ \gcd(a,q)=1}}\chi(ma+n\overline{a}),\end{equation}
which is the analogue of the Kloosterman sums. Here $\mathcal{I}$ is a subinterval of $[x+1,x+q]$ for some integer $x.$

For $\mathcal{I}=[x+1,x+q]$, the sum (\ref{eq:3}) is said to be \textit{complete}, in which case we write $\Lambda_\chi(m,n,\mathcal{I};q)=\Lambda_\chi(m,n;q)$ for short, and for $\mathcal{I}$ being the proper subset of $[x+1,x+q]$, the sum (\ref{eq:3}) is said to be \textit{partial} (or \textit{incomplete}).
Clearly, the complete sum $\Lambda_\chi(m,n;q)$ enjoy the multiplicity property
\[\Lambda_\chi(m,n;q_1q_2)=\Lambda_{\chi_1}(m,n;q_1)\Lambda_{\chi_2}(m,n;q_2)\]for $\gcd(q_1,q_2)=1,\chi_1\bmod q_1$ and $\chi_2\bmod q_2$ with $\chi_1\chi_2=\chi$.

For convenience in the following paragraphs, we introduce the definition on \emph{completely} even characters.
\begin{definition}Let $q$ be a positive integer and $\chi$ a Dirichlet character mod $q$. Suppose $q=\prod_{1\leqslant i\leqslant s}p_i^{\alpha_i}$ is the canonical decomposition of $q$ and $\chi=\prod_{1\leqslant i\leqslant s}\chi_i$ with $\chi_i\bmod p_i^{\alpha_i}$ for each $i$. If $\chi_i(-1)=1$ for each $i$, then we shall call $\chi$ a completely even character mod $q$.
\end{definition}
If $\chi(-1)=-1,$ we have \[\Lambda_\chi(m,n;q)=\sideset{}{^*}\sum_{a\bmod q}\chi(ma+n\overline{a})=\sideset{}{^*}\sum_{a\bmod q}\chi(-ma-n\overline{a})
=-\Lambda_\chi(m,n;q),\]hence $\Lambda_\chi(m,n;q)=0$. From the observation and the multiplicity property, we find that $\Lambda_\chi(m,n;q)$ must vanish unless
$\chi$ is a completely even character mod $q$, which we always assume during the investigation of $\Lambda_\chi(m,n;q)$.

For the individual sum, J. Yang and Z. Y. Zheng \cite{YZ} obtained the following upper bound for the complete sum
\begin{equation}\label{eq:4}\max_{m,n}|\Lambda_\chi(m,n;q)|\leqslant q^{1/2}2^{\omega(q)}\end{equation}
for an arbitrary integer $q\geqslant3$ and a nontrivial character $\chi\bmod q$. By Fourier technique or completing method, and the estimates for mixed exponential sums due to T. Cochrane and Z. Y. Zheng \cite{CZ}, we can deduce that
\begin{equation}\label{eq:5}\max_{m,n}|\Lambda_\chi(m,n,\mathcal{I};q)|\ll q^{1/2+\varepsilon}.\end{equation}

In order to illustrate the optimization of the estimate (\ref{eq:4}), we consider the second moment of the complete sums
\[\mathcal{K}(\chi,q)=\mathop{\sum\sum}_{m,n\bmod q}|\Lambda_\chi(m,n;q)|^2.\]
By the standard methods, we will prove in Section \ref{sec:2} an calculation formula for $\mathcal{K}(\chi,q)$. To be precise, we can state that
\begin{theorem}\label{thm:1} Let $q\geqslant3$ be a positive integer and $\chi$ a completely even primitive character mod $q$. Then we have
\[\mathcal{K}(\chi,q)=q\varphi^2(q)\displaystyle\sum_{\substack{y\bmod q\\y^2\equiv1(\bmod q)}}\chi(y).\]

Moreover, $\mathcal{K}(\chi,q)$ vanishes if $8|q$, and for $8\nmid q,$ we have
\[q\varphi^2(q)\leqslant\mathcal{K}(\chi,q)\leqslant q\varphi^2(q)2^{\omega(q)}.\]
\end{theorem}

This shows that $\mathcal{K}(\chi,q)$ is of the order of magnitude $q^{3\pm\varepsilon}$ for any completely even primitive character mod $q$ under certain assumptions, which yields the estimate (\ref{eq:4}) is rather sharp. The upper bound saves a factor $2^{\omega(q)}$ compared with the trivial conclusion from (\ref{eq:4}).

Following a similar idea to the proof of Theorem \ref{thm:1} in Section \ref{sec:2}, we can also conclude the following general result.
\begin{theorem}\label{thm:2} If $\{\lambda_a\}$ is bounded, then we have
\[\mathop{\sum\sum}_{m,n\bmod q}\left|\ \sideset{}{^*}\sum_{a\bmod q}\lambda_a\chi(ma+n\overline{a})\right|^{2}\ll q\varphi^2(q)2^{\omega(q)}.\]
\end{theorem}

In some applications, one may have to deal with the bilinear form
\[\mathcal{B}_\chi(A,M,N;\alpha,\beta,q)=\sum_{m\sim M}\sum_{n\sim N}\sideset{}{^*}\sum_{a\sim A}\alpha_m\beta_n\chi(ma+n\overline{a}),\]
where $\alpha=(\alpha_m)$ and $\beta=(\beta_n)$ are complex coefficients. From the Cauchy-Schwarz inequality and (\ref{eq:5}), one easily obtains \begin{align}|\mathcal{B}_\chi(A,M,N;\alpha,\beta,q)|&\leqslant\|\alpha\|\|\beta\|\bigg(\sum_{m\sim M}\sum_{n\sim N}\bigg|\sideset{}{^*}\sum_{a\sim A}\chi(ma+n\overline{a})\bigg|^2\bigg)^{1/2}\nonumber\\
&\ll\|\alpha\|\|\beta\|(MNq)^{1/2}q^\varepsilon,\label{eq:6}\end{align}where $\|\cdot\|$ denotes the $\ell_2$-norm. In Section \ref{sec:3}, we shall prove a much better estimate for the bilinear form.
\begin{theorem}\label{thm:3}For an arbitrary integer $q\geqslant3$ and primitive character $\chi\bmod q,$ we have
\[\mathcal{B}_\chi(A,M,N;\alpha,\beta,q)\ll\|\alpha\|\|\beta\|N^{1/2}q^{3/4}\tau^{5/2}(q)\log^2q,\] provided that each prime factor of $q$ is $\gg\log^\gamma N$ for a given $\gamma>1.$ Here the implied constant depends only on $\gamma.$
\end{theorem}

\noindent \textbf{Notation.} The following notation will be used throughout this paper.
\begin{itemize}
\item $e(x)=\exp(2\pi ix)=e^{2\pi ix};$
\item $\gcd(a,b)$ denotes the greatest common divisor of $a$ and $b$;
\item $f = O(g)$ ($f\ll g)$ means $|f|\leqslant cg$ for some unspecified positive constant $c$;
\item $x\sim X$ means $X<x\leqslant2X$;
\item $\varphi(q)$ and $\tau(q)$ denote the Euler function and divisor function, respectively; $\omega(q)$ denotes the number of distinct prime divisors of $q$ and $P(q)$ denotes the greatest prime divisor of $q$;
\item $\|\cdot\|$ denotes the $\ell_2$-norm;
\item $L^1(\mathbb{R})$ denotes the set consisting of all the complex-valued Lebesgue integrable functions over $\mathbb{R}$;
\item $\displaystyle\sideset{}{^*}\sum_{a\bmod q}$ denotes the sum restricted to a reduced residue system mod $q$;
\end{itemize}

Throughout this paper, the summations and integrals without the limitation of the range stand for the ones
over the whole range, from $-\infty$ to $+\infty$. To be precise,
\[\sum_{n}a_n:=\sum_{n\in\mathbb{Z}}a_n,\ \ \ \text{and}\ \ \int f(x)\ud x:=\int_{\mathbb{R}}f(x)\ud x.\]

\vskip 30pt

\section{Moments of the complete sums}\label{sec:2}

In this section, we shall deal with the second mean value of the complete sums
\[\mathcal{K}(\chi,q)=\mathop{\sum\sum}_{m,n\bmod q}\left|\ \sideset{}{^*}\sum_{a\bmod q}\chi(ma+n\overline{a})\right|^2\]under the assumption that
$\chi$ is a completely even primitive character mod $q$.

First, we shall list some basic properties of Gauss sum
\[G(n,\chi)=\sum_{a\bmod q}\chi(a)e\left(\frac{na}{q}\right).\]
\begin{lemma}\label{lm:1}For an arbitrary character $\chi\bmod q$, we have \[\overline{G(1,\chi)}=\chi(-1)G(1,\overline{\chi}).\] If $\chi$ is primitive, then
\[G(n,\chi)=\overline{\chi}(n)G(1,\chi)\]and\[|G(1,\chi)|=q^{1/2}.\]
\end{lemma}

Since we have
\[\sideset{}{^*}\sum_{a\bmod q}\chi(ma+n\overline{a})=\sideset{}{^*}\sum_{c\bmod q}\chi(c)\sideset{}{^*}\sum_{\substack{a\bmod q\\ma+n\overline{a}\equiv c(\bmod q)}}1,\]
we obtain
\begin{align}\mathcal{K}(\chi,q)&=\mathop{\sideset{}{^*}\sum\sideset{}{^*}\sum}_{c,d\bmod q}\chi(c)\overline{\chi}(d)\mathop{\sideset{}{^*}\sum\sideset{}{^*}\sum}_{a,b\bmod q}\mathop{\sum\sum}_{\substack{m,n\bmod q\\ma+n\overline{a}\equiv c(\bmod q)\\mb+n\overline{b}\equiv d(\bmod q)}}1.\label{eq:7}\end{align}

Now the system of congruences
\[\begin{cases}ax+\overline{a}y&\equiv c\pmod q\\ bx+\overline{b}y&\equiv d\pmod q
\end{cases}\] is solvable in integers $x$ and $y$, if and only if\[\gcd(a^2-b^2,q)|{(ad-bc)},\] in which case the number of solutions
$(x,y)$ is $\gcd(a^2-b^2,q).$ With this observation in mind, we obtain from (\ref{eq:7}) that
\[\begin{split}
\mathcal{K}(\chi,q)&=\mathop{\sideset{}{^*}\sum\sideset{}{^*}\sum}_{c,d\bmod q}\chi(c)\overline{\chi}(d)\mathop{\sideset{}{^*}\sum\sideset{}{^*}\sum}_{\substack{a,b\bmod q\\ \gcd(a^2-b^2,q)|(ad-bc)}}\gcd(a^2-b^2,q)\\
&=\sum_{\ell|q}\ell\mathop{\sideset{}{^*}\sum\sideset{}{^*}\sum}_{c,d\bmod q}\chi(c)\overline{\chi}(d)\mathop{\sideset{}{^*}\sum\sideset{}{^*}\sum}_{\substack{a,b\bmod q\\ a^2\equiv b^2(\bmod\ell)\\ ad\equiv bc(\bmod\ell)\\ \gcd(\frac{a^2-b^2}{\ell},\frac{q}{\ell})=1}}1\\
&=\sum_{\ell|q}\ell\sum_{y\bmod\ell}\mathop{\sideset{}{^*}\sum\sideset{}{^*}\sum}_{\substack{a,b\bmod q\\ a^2\equiv b^2(\bmod\ell)\\ a\equiv by(\bmod\ell)\\ \gcd(\frac{a^2-b^2}{\ell},\frac{q}{\ell})=1}}\mathcal{C}(y;\ell,q),\end{split}\]
where\[\mathcal{C}(y;\ell,q):=\mathop{\sideset{}{^*}\sum\sideset{}{^*}\sum}_{\substack{c,d\bmod q\\ c\equiv dy(\bmod \ell)}}\chi(c)\overline{\chi}(d).\]

The next goal is to calculate the character sum $\mathcal{C}(y;\ell,q)$ for each $\ell|q$ and $y\bmod q$. For each divisor $\ell$ of $q,$
\begin{align*}\mathcal{C}(y;\ell,q)
&=\frac{1}{\ell}\sum_{k\bmod \ell}\ \sideset{}{^*}\sum_{c\bmod q}\chi(c)e\left(\frac{kc}{\ell}\right)\sideset{}{^*}\sum_{d\bmod q}\overline{\chi}(d)e\left(\frac{-kdy}{\ell}\right)\\
&=\frac{1}{\ell}\sum_{k\bmod \ell}\ \sideset{}{^*}\sum_{c\bmod q}\chi(c)e\left(\frac{kcq/\ell}{q}\right)\sideset{}{^*}\sum_{d\bmod q}\overline{\chi}(d)e\left(\frac{-kdyq/\ell}{q}\right).\end{align*}
Since $\chi$ is primitive mod $q$, thus from Lemma \ref{lm:1}, we obtain that
\[\mathcal{C}(y;\ell,q)=\frac{q}{\ell}\chi^0(q/\ell)\chi(y)\sum_{\substack{k\bmod \ell\\ \gcd(k,q)=1}}1,\]
where $\chi^0$ is the trivial character mod $q$. Therefore, $\mathcal{C}(y;\ell,q)$ must vanish unless $\ell=q$ and $\gcd(y,q)=1$, in which case it equals to $\chi(y)\varphi(q)$.

Hence if $\chi$ is primitive mod $q$ and $\chi(-1)=1$, we have
\begin{align*}\mathcal{K}(\chi,q)&=q\varphi(q)\sum_{y\bmod q}\chi(y)\mathop{\sideset{}{^*}\sum\sideset{}{^*}\sum}_{\substack{a,b\bmod q\\ a^2\equiv b^2(\bmod q)\\ a\equiv by(\bmod q)}}1\\
&=q\varphi^2(q)\sum_{\substack{y\bmod q\\y^2\equiv1(\bmod q)}}\chi(y).\end{align*}

We have the multiplicative property \[\sum_{\substack{y\bmod q\\y^2\equiv1(\bmod q)}}\chi(y)=\sum_{\substack{x\bmod q_1\\x^2\equiv1(\bmod q_1)}}\chi_1(x)\sum_{\substack{z\bmod q_2\\z^2\equiv1(\bmod q_2)}}\chi_2(z)\] for $q=q_1q_2,\gcd(q_1,q_2)=1,\chi\bmod q_1,\chi_2\bmod q_2$ and $\chi_1\chi_2=\chi.$ Now the problem reduces to the case of prime power moduli.

For $q=p^s,s\geqslant1,p\geqslant3,$ the equation $y^2\equiv1(\bmod q)$ is equivalent to $y\equiv\pm1(\bmod q),$ in which case
\[\sum_{\substack{y\bmod q\\y^2\equiv1(\bmod q)}}\chi(y)=2.\] And for $q=2^s,s\geqslant3,$ the equation $y^2\equiv1(\bmod q)$ is equivalent to $y\equiv\pm1$ or $\pm1+2^{s-1}(\bmod q),$ in which case
\[\sum_{\substack{y\bmod q\\y^2\equiv1(\bmod q)}}\chi(y)=2+2\chi(1+2^{s-1}).\] From the definition of primitive characters we find that
$\chi(1+2^{s-1})\neq1,$ thus $\chi(1+2^{s-1})=-1$, which yields \[\sum_{\substack{y\bmod q\\y^2\equiv1(\bmod q)}}\chi(y)=0\] for $q=2^s,s\geqslant3$. And for $q=2$ or $4$, it is easy to check that the character sum is equal to 1 or 2 correspondingly.

To sum up, we have obtained
\[\sum_{\substack{y\bmod q\\y^2\equiv1(\bmod q)}}\chi(y)=\begin{cases}
1,& q=2,\\
2,& q=4\ \text{or}\ p^s, s\geqslant1,p\geqslant3,\\
0,& q=2^s, s\geqslant3.
\end{cases}\]
Hence we have $\mathcal{K}(\chi,q)=0$ if $8|q$, and for $8\nmid q,$ we have
\[q\varphi^2(q)\leqslant\mathcal{K}(\chi,q)\leqslant q\varphi^2(q)2^{\omega(q)}.\]

This completes the proof of Theorem \ref{thm:1}.

\vskip 30pt

\section{Estimate for the bilinear form concerning the incomplete sums}\label{sec:3}

This section focuses on the upper bound of the bilinear form
\[\mathcal{B}_\chi(A,M,N;\alpha,\beta,q)=\sum_{m\sim M}\sum_{n\sim N}\sideset{}{^*}\sum_{a\sim A}\alpha_m\beta_n\chi(ma+n\overline{a}).\] We always assume that $\chi$ is primitive mod $q$ in this section.

 \subsection{Some notations and preliminary results }We first introduce the Fourier transform and Poisson summation formula. Suppose $f\in L^1(\mathbb{R}),$  the Fourier transform of $f$ is defined by
\[\widehat{f}(\lambda)=\int f(x)e(-\lambda x)\ud x.\] If $f,\widehat{f}\in L^1(\mathbb{R}),$ and have bounded variance, then we have the following Poisson summation formula
\begin{equation}\label{eq:8}\sum_nf(n)=\sum_m\widehat{f}(m).\end{equation}In fact, we shall use an extended version of
(\ref{eq:8}).

\begin{lemma}[See \cite{IK}, p.70]\label{lm:2}
Let $f,\widehat{f}\in L^1(\mathbb{R}),$ and have bounded variance, $v>0$ and $u$ be fixed real numbers. Then
\[\sum_nf(vn+u)=\frac{1}{v}\sum_m\widehat{f}\left(\frac{m}{v}\right)e\left(\frac{um}{v}\right).\]
\end{lemma}

Now we turn to prove Theorem \ref{thm:3}. Let $\phi$ be a smooth function supported in $[A,2A]$ with value 1. We can assume the derivations satisfy
\[\phi^{(\ell)}(x)\ll A^{-\ell},\ \ \ \ell\geqslant0,\] the implied constant in $\ll$ depends only on $\ell.$
We have $\widehat{\eta}(0)\ll A$, and from the alternative integration by parts, we obtain for $\lambda\neq0$ that
\begin{align*}\widehat{\phi}(\lambda)&=(-2\pi i\lambda)^{-k}\int_A^{2A}\eta^{(k)}(x)e(-\lambda x)\ud x \nonumber \\
&\ll|\lambda|^{-k}\int_A^{2A}x^{-k}\ud x \\
&\ll|\lambda|^{-k}A^{1-k}\end{align*} for each $k\geqslant1.$ Hence we can find that
\begin{align}\widehat{\phi}(\lambda)\ll A(1+|\lambda|A)^{-k}\label{eq:9}\end{align} for any $k\geqslant0.$

Now we introduce another test function $\xi$ which is smooth and supported in $[M,2M]$ with value 1.
Suppose its derivatives satisfy
\[\xi^{(\ell)}(x)\ll M^{-\ell},\ \ \ \ell\geqslant0,\] the implied constant in $\ll$ depends only on $\ell.$
Moreover, we have $\widehat{\xi}(\lambda)\ll M(1+|\lambda|M)^{-k}$ for any $k\geqslant0.$

We require an identity concerning primitive characters.
\begin{lemma}\label{lm:3}Let $q$ be a positive integer and $\chi$ a primitive character mod $q$. Then we have
\[\frac{1}{\varphi(q)}\ \sideset{}{^*}\sum_{a\bmod q}\chi(ca+b)=\begin{cases}
\chi(b),& q\mid c,\\
0,& q\nmid c.
\end{cases}\]
\end{lemma}

\proof By the M\"{o}bius inversion formula, we obtain that
\begin{align*}\frac{1}{\varphi(q)}\ \sideset{}{^*}\sum_{a\bmod q}\chi(ca+b)&=\frac{1}{\varphi(q)}\sum_{d|q}\mu(d)\sum_{a\bmod (q/d)}\chi(cda+b)\\
&=\frac{1}{\varphi(q)}\sum_{d|q}\mu(d)\sum_{\ell\bmod q}\chi(\ell)\sum_{\substack{a\bmod (q/d)\\cda+b\equiv \ell(\bmod q)}}1\\
&=\frac{1}{q\varphi(q)}\sum_{k\bmod q}e\left(\frac{kb}{q}\right)\sum_{d|q}\mu(d)\sum_{\ell\bmod q}\chi(\ell)e\left(\frac{-k\ell}{q}\right)\\
&\ \ \ \ \ \times\sum_{a\bmod (q/d)}e\left(\frac{kca}{q/d}\right)\\
&=\frac{1}{\varphi(q)}\sum_{k\bmod q}e\left(\frac{kb}{q}\right)\sum_{\substack{d|q\\ q/d|kc}}\frac{\mu(d)}{d}\sum_{\ell\bmod q}\chi(\ell)e\left(\frac{-k\ell}{q}\right).\end{align*}

Since $\chi$ is primitive mod $q$, then from Lemma \ref{lm:3}, we have
\[\begin{split}\frac{1}{\varphi(q)}\ \sideset{}{^*}\sum_{a\bmod q}\chi(ca+b)
&=\frac{\chi(-1)G(\chi)}{\varphi(q)}\sum_{\substack{d|q\\ q/d|c}}\frac{\mu(d)}{d}\sum_{k\bmod q}\overline{\chi}(k)e\left(\frac{kb}{q}\right)\\
&=\chi(b)\frac{q}{\varphi(q)}\sum_{\substack{d|q\\ q/d|c}}\frac{\mu(d)}{d}.
\end{split}\]

In fact, \[W:=\sum_{\substack{d|q\\ q/d|c}}\frac{\mu(d)}{d}=\prod_{p^\alpha\|q}\left(f_c(p^\alpha)-\frac{f_c(p^{\alpha-1})}{p}\right),\]where
\[f_c(p^k)=\begin{cases}
1,& p^k\mid c,\\
0,& p^k\nmid c,
\end{cases}\]for each prime $p$ and positive integer $k$.

Hence for each $p$ with $p^\alpha\|q$, $W$ must vanish if $p^\alpha\nmid q$. Hence the sum in question must vanish if $q\nmid c$. The case for $q\mid c$ is also valid
on observing that $\varphi(q)=\sum_{d\ell=q}d\mu(\ell)$ as required.\endproof

We also require an estimate for complete exponential sums for polynomials.
\begin{lemma}\label{lm:4}Let $q$ be a positive integer and $f(x)=ax^2+bx\in\mathbb{Z}[x]$. Then we have
\[\sideset{}{^*}\sum_{x\bmod q}e\bigg(\frac{f(x)}{q}\bigg)\ll q^{1/2}\gcd(a,q)^{1/2}2^{\omega(q)}.\]

\end{lemma}
\proof According to the estimate due to N. M. Korobov \cite{K},
\[\sum_{x\bmod q}e\bigg(\frac{ax^2+bx}{q}\bigg)\ll q^{1/2}\gcd(a,q)^{1/2},\]
we derive from the M\"{o}bius inversion formula that
\begin{align*}\sideset{}{^*}\sum_{x\bmod q}e\bigg(\frac{ax^2+bx}{q}\bigg)&=\sum_{d|q}\mu(d)\sum_{\substack{x\bmod q\\ d|x}}e\bigg(\frac{ax^2+bx}{q}\bigg)\\
&=\sum_{d|q}\mu(d)\sum_{x\bmod q/d}e\bigg(\frac{adx^2+bx}{q/d}\bigg)\\
&\ll q^{1/2}\sum_{d|q}\mu^2(d)d^{-1/2}\gcd(ad,q/d)^{1/2}\\
&\leqslant q^{1/2}\gcd(a,q)^{1/2}\sum_{d|q}\mu^2(d)\\
&= q^{1/2}\gcd(a,q)^{1/2}2^{\omega(q)}.
\end{align*}This completes the proof of the lemma.
\endproof

\subsection{Completing the summation over $m$}

Now we begin to estimate the bilinear form.
From Lemma \ref{lm:2}, we have
\begin{align*}\sideset{}{^*}\sum_{a\sim A}\chi(ma+n\overline{a})&=\sideset{}{^*}\sum_{a}\phi(a)\chi(ma+n\overline{a})\\
&=\sideset{}{^*}\sum_{a\bmod q}\chi(ma+n\overline{a})\sum_{b}\phi(a+bq) \\
&=\frac{1}{q}\sum_\nu\widehat{\phi}\left(\frac{\nu}{q}\right)\sideset{}{^*}\sum_{a\bmod q}\chi(ma+n\overline{a})e\left(\frac{a\nu}{q}\right).
\end{align*}
Hence we have
\begin{equation}\label{eq:10}\mathcal{B}_\chi(A,M,N;\alpha,\beta,q)=\frac{1}{q}\sum_\nu\widehat{\phi}\left(\frac{\nu}{q}\right)\mathcal{B}^*_\chi(M,N;\alpha,\beta,\nu,q),
\end{equation}with
\[\mathcal{B}^*_\chi(M,N;\alpha,\beta,\nu,q)=\sum_{m\sim M}\sum_{n\sim N}\alpha_m\beta_n\sideset{}{^*}\sum_{a\bmod q}\chi(ma+n\overline{a})e\left(\frac{a\nu}{q}\right)=:\mathcal{B}^*_\chi.\]

From the Cauchy-Schwarz inequality and Lemma \ref{lm:3}, we have
\begin{align*}|\mathcal{B}^*_\chi|^2&\leqslant\|\alpha\|^2\sum_{m\sim M}\left|\sum_{n\sim N}\beta_n\sideset{}{^*}\sum_{a\bmod q}\chi(ma+n\overline{a})e\left(\frac{a\nu}{q}\right)\right|^2\end{align*}

Recalling the definition of $\xi$, we have
\begin{align}|\mathcal{B}^*_\chi|^2&\leqslant\|\alpha\|^2\sum_{m}\xi(m)\left|\sum_{n\sim N}\beta_n\sideset{}{^*}\sum_{a\bmod q}\chi(ma+n\overline{a})e\left(\frac{a\nu}{q}\right)\right|^2\nonumber\\
&=\frac{\|\alpha\|^2}{q}\sum_k\widehat{\xi}\left(\frac{k}{q}\right)\sum_{m\bmod q}\left|\sum_{n\sim N}\beta_n\sideset{}{^*}\sum_{a\bmod q}\chi(ma+n\overline{a})e\left(\frac{a\nu}{q}\right)\right|^2e\left(\frac{km}{q}\right)\nonumber
\end{align}

From the upper bound (\ref{eq:9}), we have
\[\sum_k\left|\widehat{\xi}\left(\frac{k}{q}\right)\right|\ll q\log q,\]
thus we can get
\begin{align}|\mathcal{B}^*_\chi|^2&\ll\|\alpha\|^2\log q\sum_{m\bmod q}\left|\sum_{n\sim N}\beta_n\sideset{}{^*}\sum_{a\bmod q}\chi(ma+n\overline{a})e\left(\frac{a\nu}{q}\right)\right|^2\nonumber\\
&=\|\alpha\|^2\log q\mathop{\sum\sum}_{n_1,n_2\sim N}\beta_{n_1}\overline{\beta}_{n_2}\mathop{\sideset{}{^*}\sum\sideset{}{^*}\sum}_{a,b\bmod q}e\left(\frac{\nu(a-b)}{q}\right)\nonumber\\
&\ \ \ \ \ \ \ \ \ \ \times\sum_{m\bmod q}\chi(ma+n_1\overline{a})\overline{\chi}(mb+n_2\overline{b}).\label{eq:11}
\end{align}

Changing variables iteratively, we obtain
\[\begin{split}\sum_{m\bmod q}\chi(ma+n_1\overline{a})\overline{\chi}(mb+n_2\overline{b})&=\sum_{m\bmod q}\chi(m+n_1\overline{a})\overline{\chi}(m\overline{a}b+n_2\overline{b})\\
&=\sum_{m\bmod q}\chi(m)\overline{\chi}(m\overline{a}b+n_2\overline{b}-n_1\overline{a}^2b)\\
&=\sideset{}{^*}\sum_{m\bmod q}\overline{\chi}(\overline{a}b+(n_2\overline{b}-n_1\overline{a}^2b)m).
\end{split}\]From Lemma \ref{lm:3}, we obtain that the summation above vanishes unless $n_2a^2\equiv n_1b^2\pmod q$, in which case the sum can be reduced to $\varphi(q)\chi(a\overline{b}).$ Hence the right-hand side of (\ref{eq:11}) is just
\begin{align}\|\alpha\|^2\varphi(q)\log q\mathop{\sum\sum}_{n_1,n_2\sim N}&\beta_{n_1}\overline{\beta}_{n_2}\mathop{\sideset{}{^*}\sum\sideset{}{^*}\sum}_{\substack{a,b\bmod q\\ n_2a^2\equiv n_1b^2(\bmod q)}}\chi(a\overline{b})e\left(\frac{\nu(a-b)}{q}\right)\nonumber\\
&\leqslant\|\alpha\|^2\|\beta\|^2\varphi(q)\log q\sum_{L\ll\log q}\mathcal{J}(L,\beta;\nu,q),\label{eq:12}\end{align}
where
\[\mathcal{J}(L,\beta;\nu,q)=\sum_{\ell\sim e^L}\left|\sum_{n\sim N}{\beta}_{n}\sideset{}{^*}\sum_{\substack{a\bmod q\\ na^2\equiv\ell(\bmod q)}}\chi(a)e\left(\frac{\nu a}{q}\right)\right|^2.\]

\subsection{Estimate for $\mathcal{J}(L,\beta;\nu,q)$} We shall estimate $\mathcal{J}(L,\beta;\nu,q)$ beginning with the ideas from the sieve technique. For each $\ell$, we have
\[\begin{split}\sum_{n\sim N}{\beta}_{n}\sideset{}{^*}\sum_{\substack{a\bmod q\\ na^2\equiv\ell(\bmod q)}}\chi(a)e\left(\frac{\nu a}{q}\right)
=\sum_{\substack{n\sim N\\ \gcd(n,r)=1}}&\beta_{n}\sideset{}{^*}\sum_{\substack{a\bmod q\\ na^2\equiv\ell(\bmod q)}}\chi(a)e\left(\frac{\nu a}{q}\right)\\
&+\sum_{n\sim N}\beta'_{n}\sideset{}{^*}\sum_{\substack{a\bmod q\\ na^2\equiv\ell(\bmod q)}}\chi(a)e\left(\frac{\nu a}{q}\right),\end{split}\]
where $r$ is the product of primes,
\[r=\prod_{\log^\gamma N<p\leqslant R}p,\ \ \ \ R>2N.\]Here $\gamma>1$ is a given constant, $R\in\mathcal{R}$ is a parameter to be chosen later, $\mathcal{R}$ is a finite set such that
\[\sum_{\substack{R\in\mathcal{R},R>2N\\p|r}}1\ll\frac{1}{p}\sum_{R\in\mathcal{R},R>2N}1\]holds for any prime $p$ with $p>\log^\gamma N$, and
\[\beta'_{n}=\begin{cases}
\beta_{n},& \text{if}\ \gcd(n,r)>1,\\
0,& \text{otherwise}.
\end{cases}
\]

From the definition of $r$ we see that $\gcd(n,r)=1$ is equivalent to $P(n)\leqslant\log^\gamma N,$ where $P(n)$ denotes the greatest prime factor of $n$. Hence we obtain from the inequality $|x+y|^2\leqslant2(|x|^2+|y|^2)$ that
\[\begin{split}\left|\sum_{n\sim N}{\beta}_{n}\sideset{}{^*}\sum_{\substack{a\bmod q\\ na^2\equiv\ell(\bmod q)}}\chi(a)e\left(\frac{\nu a}{q}\right)\right|^2
&\leqslant2\left|\sum_{\substack{n\sim N\\ P(n)\leqslant\log^\gamma N}}\beta_{n}\sideset{}{^*}\sum_{\substack{a\bmod q\\ na^2\equiv\ell(\bmod q)}}\chi(a)e\left(\frac{\nu a}{q}\right)\right|^2\\
&\ \ \ \ \ \ \ +2\left|\sum_{n\sim N}\beta'_{n}\sideset{}{^*}\sum_{\substack{a\bmod q\\ na^2\equiv\ell(\bmod q)}}\chi(a)e\left(\frac{\nu a}{q}\right)\right|^2.\end{split}\]
Summing over $\ell$ we arrive at
\begin{equation}\label{eq:13}\mathcal{J}(L,\beta;\nu,q)\leqslant2\mathcal{J}(R,L,\beta;\nu,q)+2\mathcal{J}(L,\beta';\nu,q),\end{equation}
where $\beta'=(\beta'_{n}),$ and \[\mathcal{J}(R,L,\beta;\nu,q)=\sum_{\ell\sim e^L}\left|\sum_{\substack{n\sim N\\ P(n)\leqslant\log^\gamma N}}\beta_{n}\sideset{}{^*}\sum_{\substack{a\bmod q\\ na^2\equiv\ell(\bmod q)}}\chi(a)e\left(\frac{\nu a}{q}\right)\right|^2.\]

Suppose $\eta(\ell)$ is a smooth function supported on $(e^L/2,e^L]$ such that $\eta^{(i)}(\ell)\ll e^{-iL}$ for all $i\geqslant0$. Then
\[\begin{split}\mathcal{J}(R,L,\beta;\nu,q)&\leqslant\sum_{\ell}\eta(\ell)\left|\sum_{\substack{n\sim N\\ P(n)\leqslant\log^\gamma N}}\beta_{n}\sideset{}{^*}\sum_{\substack{a\bmod q\\ na^2\equiv\ell(\bmod q)}}\chi(a)e\left(\frac{\nu a}{q}\right)\right|^2.\end{split}\]

After squaring out and changing the order of summation, we
obtain from the Poisson summation formula that
\[\begin{split}\mathcal{J}(R,L,\beta;\nu,q)&\leqslant\mathop{\sum\sum}_{\substack{n_1,n_2\sim N\\ P(n_1n_2)\leqslant\log^\gamma N}}\beta_{n_1}\overline{\beta}_{n_2}\mathop{\sideset{}{^*}\sum\sideset{}{^*}\sum}_{\substack{a,b\bmod q\\ n_1a^2\equiv n_2b^2(\bmod q)}}\chi(a\overline{b})e\left(\frac{\nu(a-b)}{q}\right)\sum_{\ell\equiv n_1a^2(\bmod q)}\eta(\ell)\\
&=\frac{1}{q}\mathop{\sum\sum}_{\substack{n_1,n_2\sim N\\ P(n_1n_2)\leqslant\log^\gamma N}}\beta_{n_1}\overline{\beta}_{n_2}\mathop{\sideset{}{^*}\sum\sideset{}{^*}\sum}_{\substack{a,b\bmod q\\ n_1a^2\equiv n_2b^2(\bmod q)}}\chi(a\overline{b})e\left(\frac{\nu(a-b)}{q}\right)\sum_{\ell}\widehat{\eta}\left(\frac{\ell}{q}\right)e\left(\frac{n_1a^2\ell}{q}\right)\\
&=\frac{1}{q}\mathop{\sum\sum}_{\substack{n_1,n_2\sim N\\ P(n_1n_2)\leqslant\log^\gamma N}}\beta_{n_1}\overline{\beta}_{n_2}\sum_{\ell}\widehat{\eta}\left(\frac{\ell}{q}\right)\sum_{\substack{y\bmod q\\n_1y^2\equiv n_2(\bmod q)}}\chi(y)\sideset{}{^*}\sum_{a\bmod q}e\left(\frac{\ell n_1a^2+\nu a(1-\overline{y})}{q}\right).\end{split}\]

Since $q$ has no small prime divisors, that is each prime factor of $q$ is $\gg\log^\gamma N,$ then we have $\gcd(n_1n_2,q)=1$ for $P(n_1n_2)\leqslant\log^\gamma N.$
From Lemma \ref{lm:4} and the estimate
\[\sum_{\ell\geqslant1}\left|\widehat{\eta}\left(\frac{\ell}{q}\right)\right|\gcd(\ell,q)^{1/2}\ll q\tau(q)\log q,\]
we find that
\begin{align}\mathcal{J}(R,L,\beta;\nu,q)&\ll e^L\mathop{\sum\sum}_{\substack{n_1,n_2\sim N\\\gcd(n_1n_2,q)=1}}|\beta_{n_1}||\beta_{n_2}|\sideset{}{^*}\sum_{\substack{y\bmod q\\n_1y^2\equiv n_2(\bmod q)}}1\nonumber\\
&\ \ \ \ \ \ \ \ \ +q^{-1/2}4^{\omega(q)}\mathop{\sum\sum}_{\substack{n_1,n_2\sim N\nonumber\\ \gcd(n_1n_2,q)=1}}|\beta_{n_1}||\beta_{n_2}|\sum_{\ell\geqslant1}\left|\widehat{\eta}\left(\frac{\ell}{q}\right)\right|\gcd(\ell,q)^{1/2}\nonumber\\
&\ll \|\beta\|^2Nq^{1/2}\tau^3(q)\log q.\label{eq:14}\end{align}

Let $\mathcal{V}(e^L,N)$ denote the norm of the linear operator given by the matrix with the $(\ell,n)$-th entry
\[\sideset{}{^*}\sum_{\substack{a\bmod q\\ na^2\equiv\ell(\bmod q)}}\chi(a)e\left(\frac{\nu a}{q}\right)\]with $\ell\sim e^L$ and $n\sim N$. Thus
\[\mathcal{J}(R,L,\beta;\nu,q)\leqslant\|\beta\|^2\mathcal{V}(e^L,N)\]for any complex numbers $(\beta_n)$. This together with (\ref{eq:13}) and (\ref{eq:14}) yields
\begin{equation}\label{eq:15}\mathcal{J}(L,\beta;\nu,q)\ll\|\beta\|^2Nq^{1/2}\tau^3(q)\log q+\|\beta'\|^2\mathcal{V}(e^L,N).\end{equation}

By the definition of $\beta'$, we have
\begin{align*}\sum_{R\in\mathcal{R},R>2N}\|\beta'\|^2&=\sum_{R\in\mathcal{R},R>2N}\sum_{\substack{n\sim N\\ \gcd(n,r)>1}}|\beta_n|^2=\sum_{n\sim N}|\beta_n|^2\sum_{\substack{R\in\mathcal{R},R>2N\\ \gcd(n,r)>1}}1\\
&=\sum_{n\sim N}|\beta_n|^2\sum_{\substack{p|n\\ p\geqslant\log^\gamma N}}\sum_{\substack{R\in\mathcal{R},R>2N\\ p|r}}1\\
&\ll\sum_{n\sim N}|\beta_n|^2\sum_{\substack{p|n\\ p\geqslant\log^\gamma N}}\frac{1}{p}\sum_{R\in\mathcal{R},R>2N}1\\
&\leqslant\bigg(\sum_{R\in\mathcal{R},R>2N}1\bigg)\sum_{n\sim N}|\beta_n|^2\frac{\omega(n)}{\log^\gamma N}\\
&\leqslant(\log N)^{1-\gamma}\bigg(\sum_{R\in\mathcal{R},R>2N}1\bigg)\|\beta\|^2.\end{align*} Hence we can always choose $R\in\mathcal{R},R>2N$ such that
\[\|\beta'\|^2\ll \|\beta\|^2(\log N)^{1-\gamma},\] from which and (\ref{eq:15}), we obtain that
\[\mathcal{J}(L,\beta;\nu,q)\ll\|\beta\|^2(Nq^{1/2}\tau^2(q)\log q+\mathcal{V}(e^L,N)(\log N)^{1-\gamma})\]holds for any complex numbers $(\beta_n)$. In other words,
\[\mathcal{V}(e^L,N)\ll Nq^{1/2}\tau^3(q)\log q+\mathcal{V}(e^L,N)(\log N)^{1-\gamma},\]
which gives
\[\mathcal{V}(e^L,N)\ll Nq^{1/2}\tau^3(q)\log q.\]

Hence we conclude the estimate
\begin{equation}\label{eq:16}\mathcal{J}(L,\beta;\nu,q)\ll \|\beta\|^2Nq^{1/2}\tau^3(q)\log q.\end{equation}

\subsection{Conclusion}

Combining (\ref{eq:12}) and (\ref{eq:16}) we obtain that
\[\mathcal{B}_\chi^*\ll\|\alpha\|\|\beta\|N^{1/2}q^{3/4}\tau^{3/2}(q)\log q,\]
which together with (\ref{eq:10}) yields
\[\mathcal{B}_\chi(A,M,N;\alpha,\beta,q)\ll\|\alpha\|\|\beta\|N^{1/2}q^{3/4}\tau^{5/2}(q)\log^2q.\]
This completes the proof of Theorem \ref{thm:3}.

\vskip 30pt

\section{Final remarks}\label{sec:4}
Regarding the bilinear form, if each prime factor of $q$ is of the size $\gg\log^\gamma(M+N),$ then by symmetry, we have
\[\mathcal{B}_\chi(A,M,N;\alpha,\beta,q)\ll\|\alpha\|\|\beta\|\min(M^{1/2},N^{1/2})q^{3/4}\tau^{5/2}(q)\log^2q,\]which is better than (\ref{eq:6}) in the case of $M+N\gg q^{1/2+\varepsilon}$.
Also, if $\beta=(\beta_n)$ is supported on the integers $n$ with $\gcd(n,q)=1$, the proof of Theorem \ref{thm:3} becomes much simpler and it is not necessary for $q$ to have no small divisors. One finds that the estimate in Theorem \ref{thm:3} also holds for arbitrary $\alpha$ and such $\beta$.

\bigskip

\noindent \textbf{Acknowledgement.}  The last section of this paper was completed after I learned of the joint work of W. Duke, J. Friedlander and H. Iwaniec \cite{DFI}, to whom I would like to express my sincere thanks. Recently, Prof. W.P. Zhang has independently investigated $\Lambda_\chi(m,n;q)$ and found an interesting identity under certain assumptions on the moduli $q$ and the character $\chi$, I would like to thank him here for sharing his ideas during Workshop on Number Theory --- Xi'an 2011. I am also profoundly grateful to my supervisor Professor Y. Yi for the constant help, encouragement and meticulous guidance.

\bibliographystyle{plainnat}

\bigskip

\end{document}